\documentclass{gtart_h}
\def\ifplaintex{\expandafter\ifx\csname documentclass\endcsname\relax}

\def\gtp{{\mathsurround=0pt\it $\cal G\mskip-2mu$eometry \&\ 
$\cal T\!\!$opology $\cal P\!$ublications}}  % GT publications

\def\recd{{\small Received:\qua\receiveddate\ifx\reviseddate\relax
\else\qquad Revised:\qua\reviseddate\fi\par}} 

\def\lognumber#1{\def\thelognumber{#1}}
\def\volumenumber#1{\def\thevolumenumber{#1}}
\def\volumeyear#1{\def\thevolumeyear{#1}}
\def\papernumber#1{\def\thepapernumber{#1}}
\def\pagenumbers#1#2{\def\startpage{#1}\def\finishpage{#2}}
\def\published#1{\def\publishdate{#1}}

\def\received#1{\def\receiveddate{#1}}

\def\accepted#1{\def\accepteddate{#1}}

\def\asciiauthors#1{\def\theasciiauthors{#1}}
\def\asciiaddress#1{\def\theasciiaddress{#1}}
\def\asciiemail#1{\def\theasciiemail{#1}}

\long\def\asciiabstract#1{\long\def\theasciiabstract{#1}}

\let\\\par\let\thelognumber\relax\let\thevolumenumber\relax
\let\thepapernumber\relax\let\thevolumeyear\relax\let\startpage\relax
\let\finishpage\relax\let\publishdate\relax\let\receiveddate\relax
\let\reviseddate\relax\let\accepteddate\relax\let\theasciititle\relax
\let\theasciiauthors\relax\let\theasciiaddress\relax
\let\theasciiabstract\relax

\let\theasciiemail\relax

\ifplaintex
\font\logobig=cmssbx10 scaled 3836
\font\logomed=cmssbx10 scaled 2557
\else
\font\logobig=cmssbx10 scaled 4200
\font\logomed=cmssbx10 scaled 2800
\fi

\long\def\makeagttitle{   %%% start of definition of \makeagttitle
\count0=\startpage
\agt\hfill      %   Journal title (top left) 
\hbox to 45truept{\vbox to 0pt{\vglue -13truept{\logomed A\kern -.37em{\logobig 
T}\kern -.38em G}\vss}\hss}
\break
{\small Volume \thevolumenumber\ (\thevolumeyear)
\startpage--\finishpage\nl
Published: \publishdate}

\vglue .25truein

{\parskip=0pt\leftskip 0pt plus
1fil\def\\{\par\smallskip}{\Large\bf\thetitle}\par\medskip} \vglue
0.05truein

{\parskip=0pt\leftskip 0pt plus 1fil\def\\{\par}{\sc\theauthors}
\par\medskip}%
 
\vglue 0.03truein 

{\small\leftskip 25truept\rightskip 25truept{\bf Abstract}\stdspace\theabstract

{\bf AMS Classification}\stdspace\theprimaryclass
\ifx\thesecondaryclass\relax\else; \thesecondaryclass\fi\par
{\bf Keywords}\stdspace \thekeywords\par}\vglue 7truept

}   %%%% end of definition of \makeagttitle

\ifplaintex
\font\phead=cmsl9 scaled 950
\font\pnum=cmbx10 scaled 913
\font\pfoot=cmsl9 scaled 950
\headline{\vbox to 0pt{\vskip -4.5mm\line{\small\phead\ifnum
\count0=\startpage ISSN 1472-2739 (on-line) 1472-2747 (printed)
\hfill {\pnum\folio}\else\ifodd\count0\def\\{ }% 
\ifx\theshorttitle\relax\thetitle\else\theshorttitle\fi\hfill{\pnum\folio}
\else\def\\{ and }{\pnum\folio}\hfill\ifx\theshortauthors\relax\theauthors
\else\theshortauthors\fi\fi\fi}\vss}}
\footline{\vbox to 0pt{\vglue 0mm\line{\small\pfoot\ifnum\count0=\startpage
\copyright\ \gtp\hfill\else
\agt, Volume \thevolumenumber\ (\thevolumeyear)\hfill\fi}\vss}}
\else
\font=cmsl9 scaled 1050
\font\lnum=cmbx10 
\font=cmsl9 scaled 1050
\makeatletter
\def\@oddhead{{\small\lhead\ifnum\count0=\startpage ISSN 1472-2739 
(on-line) 1472-2747 (printed)\hfill {\lnum\number\count0}\else\ifodd\count0
\def\\{ }\ifx\theshorttitle\relax \thetitle \else\theshorttitle\fi\hfill
{\lnum\number\count0}\else\def\\{ and }{\lnum\number\count0}
\hfill\ifx\theshortauthors\relax 
\theauthors\else\theshortauthors\fi\fi\fi}}\def\@evenhead{\@oddhead}
\def\@oddfoot{\small\lfoot\ifnum\count0=\startpage\copyright\ \gtp\hfill\else
\agt, Volume \thevolumenumber\ (\thevolumeyear)\hfill\fi}
\def\@evenfoot{\@oddfoot}
\makeatother
\fi
\let\maketitlepage\makeagttitle
\let\makeshorttitle\maketitlepage
\let\maketitle\maketitlepage

\newwrite\gtoutfile
\long\gdef\makeheadfile{  %%% start of definition of \makeheadfile
{\def\\{, }\def\s{ }
\immediate\openout\gtoutfile head.xxx
\immediate\write\gtoutfile{Proxy-for: \ifx\theasciiauthors\relax
\theauthors\else\theasciiauthors\fi\s<\ifx\theasciiemail\relax\theemail\else\theasciiemail\fi>}
\immediate\write\gtoutfile{\noexpand\\}
\immediate\write\gtoutfile{Authors: \ifx\theasciiauthors\relax
\theauthors\else\theasciiauthors\fi}
{\def\\{ }\immediate\write\gtoutfile{Title: \ifx\theasciititle\relax
\thetitle\else\theasciititle\fi}}
\immediate\write\gtoutfile{Subj-class: GT or SG, GR etc}
\immediate\write\gtoutfile{MSC-class: \theprimaryclass\ifx\thesecondaryclass\relax\else, \thesecondaryclass\fi}
\immediate\write\gtoutfile{Journal-ref: Algebr. Geom. Topol. \thevolumenumber\s
(\thevolumeyear) \startpage-\finishpage}
\immediate\write\gtoutfile{Comments: Published by Algebraic and
Geometric Topology at}
\immediate\write\gtoutfile{\s\s\s  http://www.maths.warwick.ac.uk/agt/AGTVol\thevolumenumber/agt-\thevolumenumber-\thepapernumber.abs.html}
\immediate\write\gtoutfile{\noexpand\\}
\immediate\write\gtoutfile{}
\ifx\theasciiabstract\relax
\immediate\write\gtoutfile{\theabstract}\else
\immediate\write\gtoutfile{\theasciiabstract}\fi
\immediate\write\gtoutfile{}
\immediate\write\gtoutfile{\noexpand\\}
\immediate\write\gtoutfile{}
\immediate\closeout\gtoutfile}}  %%% end of definition of \makeheadfile

\def\maketitlepage{\makeagttitle\makeheadfile}
\let\makeshorttitle\maketitlepage
\let\maketitle\maketitlepage

\lognumber{51}
\volumenumber{5}
\volumeyear{2005}
\papernumber{51}
\pagenumbers{1315}{1324}
\received{6 July 2005} 
\accepted{2 September 2005}
\published{6 October 2005}

\usepackage{amsmath,amssymb}

\newtheorem{theorem}{Theorem}[section]
\newtheorem{proposition}[theorem]{Proposition}
\newtheorem*{mainthm}{Main theorem}
\theoremstyle{definition}
\newtheorem{remark}[theorem]{Remark}

\newcommand{\F}{\Bbb F}
\newcommand{\D}{\Delta}

\begin{document}

\title{Twisted Alexander polynomials and\\surjectivity 
of a group homomorphism}

\author{Teruaki Kitano\\Masaaki Suzuki\\Masaaki Wada}
\asciiaddress{Department of Mathematical and Computing Sciences,
Tokyo Institute of Technology\\2-12-1-W8-43 Oh-okayama, Meguro-ku,
Tokyo, 152-8552 Japan\\Graduate School of Mathematical Sciences,
The University of Tokyo\\3-8-1 Komaba, Meguro-ku, Tokyo, 153-8914 Japan
\\and\\Department of Information and Computer Sciences,
Nara Women's University\\Kita-Uoya Nishimachi, Nara, 630-8506 Japan}

\address{Department of Mathematical and Computing Sciences,
Tokyo Institute of Technology\\2-12-1-W8-43 Oh-okayama, Meguro-ku,
Tokyo, 152-8552 Japan\\\smallskip\\Graduate School of Mathematical Sciences,
The University of Tokyo\\3-8-1 Komaba, Meguro-ku, Tokyo, 153-8914 Japan
\\{\rm and}\\Department of Information and Computer Sciences,
Nara Women's University\\Kita-Uoya Nishimachi, Nara, 630-8506 Japan}

\asciiemail{kitano@is.titech.ac.jp, macky@ms.u-tokyo.ac.jp, wada@ics.nara-wu.ac.jp} 
\gtemail{\mailto{kitano@is.titech.ac.jp}, 
\mailto{macky@ms.u-tokyo.ac.jp}{\rm\qua and\qua}\mailto{wada@ics.nara-wu.ac.jp}}


\shortauthors{Teruaki Kitano, Masaaki Suzuki and Masaaki Wada}
\asciiauthors{Teruaki Kitano, Masaaki Suzuki and Masaaki Wada}

\begin{abstract}
If $\varphi \co G\rightarrow G'$ is a surjective homomorphism, 
we prove that the twisted Alexander polynomial of $G$ 
is divisible by the twisted Alexander polynomial of $G'$. 
As an application, we show non-existence of surjective homomorphism 
between certain knot groups. 
\end{abstract}
\asciiabstract{%
If phi: G-->G' is a surjective homomorphism, we prove that the twisted
Alexander polynomial of G is divisible by the twisted Alexander
polynomial of G'.  As an application, we show non-existence of
surjective homomorphism between certain knot groups.}

\maketitle

\section{Introduction}
Suppose that $G$ is a finitely presentable group 
with a surjective homomorphism to the free abelian group 
of rank $l$, eg, abelianization. 
Let $\rho \co G\rightarrow GL(n;R)$ be a linear representation. 
The twisted Alexander polynomial of $G$ associated to $\rho$ 
was introduced in \cite{Wada} 
and is defined to be a rational expression of $l$ indeterminates. 

Let $\varphi \co G\rightarrow G'$ be a surjective homomorphism. 
Each representation $\rho' \co G'\rightarrow GL(n;R)$ 
naturally induces a representation  of $G$, 
namely, $\rho=\rho'\circ\varphi$. 
In this paper we prove the following:
\begin{mainthm}
The twisted Alexander polynomial of $G$ associated to $\rho$ 
is divisible by the twisted Alexander polynomial of $G'$ 
associated to $\rho'$. 
\end{mainthm}

The corresponding fact about the Alexander polynomial is known
\cite{CF}.

We present two separate proofs of the main theorem. 
First we give a purely algebraic proof in \S \ref{alg-proof}. 
If $G$ is a knot group, 
the twisted Alexander polynomial of $G$ may be 
regarded as the Reidemeister torsion. 
In \S \ref{geo-proof}, 
we provide another proof of the main theorem 
in case when $G$ and $G'$ are knot groups, 
from the view point of the Reidemeister torsion. 

In the last section, 
we show non-existence of surjective homomorphism 
between certain knot groups, 
as an application of the main theorem. 

%%%%%%%%%%%%%%%%%%%%%%%%%%%%%%%%%%%%%%%%%%
\section{Twisted Alexander polynomial}
In this section, 
we recall briefly the definition of
the twisted Alexander polynomial. 

Let $G$ be a finitely presentable group. 
Choose and fix a presentation as follows:
\[
G =
\langle x_1,\ldots,x_u ~|~ r_1,\ldots,r_v \rangle .
\]
We denote by $\alpha \co G \to {\mathbb Z}^l$ 
a surjective homomorphism to the free abelian group 
with generators $t_1,\ldots,t_l$ 
and $\rho \co G \to GL(n;R)$ 
a linear representation, 
where $R$ is a unique factorization domain.
These maps naturally induce ring homomorphisms
$\tilde{\rho}$ and $\tilde{\alpha}$
from ${\mathbb Z}[G]$ to $M(n;R)$ and
${\mathbb Z}[{t_1}^{\pm 1},\ldots,{t_l}^{\pm 1}]$ respectively,
where
$M(n;R)$ denotes the matrix algebra of degree $n$ over $R$.
Then
$\tilde{\rho}\otimes\tilde{\alpha}$
defines a ring homomorphism
\[
{\mathbb Z}[G]\to
M\left(n;R[{t_1}^{\pm 1},\ldots,{t_l}^{\pm 1}]\right).
\]
Let
$F_u$ be the free group on
generators $x_1,\ldots,x_u$ and
$$
\Phi \co \Bbb Z[F_u]\to
M\left(n;R[{t_1}^{\pm 1},\ldots,{t_l}^{\pm 1}]\right)
$$
the composite of the surjection
${\mathbb Z}[F_u]\to{\mathbb Z}[G]$
induced by the fixed presentation and the map 
$\tilde{\rho}\otimes\tilde{\alpha} \co 
{\mathbb Z}[G]\to
M(n;R[{t_1}^{\pm 1},\ldots,{t_l}^{\pm 1}])$. 

We define the $v \times u$ matrix $M$ 
whose $(i,j)$ component is the $n \times n$ matrix
$$
\Phi\left(\frac{\partial r_i}{\partial x_j}\right)
\in M\left(n;R[{t_1}^{\pm 1},\ldots,{t_l}^{\pm 1}]\right),
$$
where
${\partial}/{\partial x}$
denotes the Fox derivation. 
This matrix $M$ is called
the Alexander matrix of
the presentation of $G$
associated to the representation $\rho$.

It is easy to see that 
there is an integer $1 \leq j \leq u$ 
such that 
$\det \Phi(x_j-1)\neq 0$. 
For such $j$, 
let us denote by $M_j$
the $v \times (u-1)$ matrix obtained from $M$
by removing the $j$-th column.
We regard $M_j$ as
an $n v \times n (u-1)$ matrix with coefficients in
$R[{t_1}^{\pm 1},\ldots,{t_l}^{\pm 1}]$.
Moreover, for an $n (u-1)$-tuple of indices
\[
I =
\left( i_1, i_2, \ldots, i_{n(u-1)} \right)
, \quad
\left( 1 \leq i_1 < i_2 < \cdots < i_{n(u-1)} \leq n v \right)
\]
we denote by $M_j^I$ the
$n (u-1) \times n (u-1)$ square matrix
consisting of the $i_k$-th row of the matrix $M_j$,
where $k=1,2,\ldots,n(u-1)$.

Then the twisted Alexander polynomial
(see \cite{Wada}) of a finitely presented group $G$
for a representation $\rho \co G \to GL(n;R)$
is defined to be a rational expression 
$$
\D_{G,\rho}(t_1,\ldots,t_l)
=\frac{\gcd_I (\det M_j^I)}{\det \Phi(x_j-1)}
$$
and
moreover is well-defined
up to a factor
$\epsilon {t_1}^{\varepsilon_1} \cdots {t_l}^{\varepsilon_l}$,
where $\epsilon \in R^{\times},\varepsilon_i \in {\mathbb Z}$. 
See \cite{Wada}, \cite{Lin}, \cite{GKM} and \cite{KL} 
for more precise definition and applications. 

\section{Main theorem and the algebraic proof}\label{alg-proof} 
In this section,
we prove the following main theorem of this paper. 

\begin{theorem}\label{main-theorem}
Let $G$ and $G'$ be finitely presentable groups and
$\alpha,\alpha'$ surjective homomorphisms
from $G,G'$ to ${\mathbb Z}^l$ respectively.
Suppose that there exists a surjective homomorphism
$\varphi \co G \to G'$ such that
$\alpha = \alpha' \circ \varphi$.
Then $\Delta_{G,\rho}$
is divisible by $\Delta_{G',\rho'}$
for any representation $\rho' \co G' \to GL(n ; R)$,
where $\rho = \rho' \circ \varphi$.
That is to say,
the quotient of $\Delta_{G,\rho}$
by $\Delta_{G',\rho'}$ is a genuine polynomial.
\end{theorem}

\begin{proof}
Choose and fix a presentation 
\[
G = \langle x_1,x_2,\ldots,x_u ~ | ~
r_1,r_2,\ldots,r_v \rangle .
\]
Since $\varphi$ is surjective, then $G'$ is generated by
$\varphi(x_1), \ldots , \varphi(x_u)$.
Namely, $G'$ can be presented as
\[
G'
=
\langle
\varphi(x_1), \varphi(x_2), \ldots , \varphi(x_u) ~ | ~
s_1,s_2,\ldots,s_{v'}
\rangle .
\]
For convenience,
we also write $x_i$ for $\varphi(x_i)$,
that is, we consider that
$G'$ is generated by $x_1,\ldots,x_u$.
By this notation, 
each relator $r_i$ is written as
\[
r_i =
\prod_k u_{k} s_{l_{i_k}}^{\varepsilon_{i_k}} u_{k}^{-1},
\quad
i = 1,2,\ldots,v, ~ 1 \leq l_{i_k} \leq v' , ~
u_k \in F_u, ~
\varepsilon_{i_k} = \pm 1 , 
\]
since $\varphi$ is a homomorphism. 
By applying the Fox derivation $\frac{\partial}{\partial x_j}$ 
and
collecting terms of $\frac{\partial s_k}{\partial x_j}$,
we get
\begin{equation}\label{r-and-s}
\varphi \left( \frac{\partial r_i}{\partial x_j} \right) =
\sum_{k=1}^{v'}
A_{i,k} \,
\frac{\partial s_k}{\partial x_j} .
\end{equation}
Here $A_{i,k}$ ($ 1 \leq i \leq v $) is
a sum of some $\varepsilon_{\bullet} \, \varphi (u_{\bullet})$,
which does not depend on $j$.
Let $M_G$ and $M_{G'}$ be
the Alexander matrices with the $u$-th column removed: 
\[
M_G
=
\left(
\begin{array}{ccc}
\tilde{\rho} \otimes \tilde{\alpha}
\left(\frac{\partial r_1}{\partial x_1}\right)
& \cdots &
\tilde{\rho} \otimes \tilde{\alpha}
\left(\frac{\partial r_1}{\partial x_{u-1}} \right)
\\
\vdots & \ddots & \vdots \\
\tilde{\rho} \otimes \tilde{\alpha}
\left(\frac{\partial r_v}{\partial x_1} \right)
& \cdots &
\tilde{\rho} \otimes \tilde{\alpha}
\left(\frac{\partial r_v}{\partial x_{u-1}} \right)
\end{array}
\right)
\]
\[
M_{G'}
=
\left(
\begin{array}{ccc}
\tilde{\rho}' \otimes \tilde{\alpha}'
\left( \frac{\partial s_1}{\partial x_1} \right)
& \cdots &
\tilde{\rho}' \otimes \tilde{\alpha}'
\left( \frac{\partial s_1}{\partial x_{u-1}} \right) \\
\vdots & \ddots & \vdots \\
\tilde{\rho}' \otimes \tilde{\alpha}'
\left( \frac{\partial s_{v'}}{\partial x_1} \right)
& \cdots &
\tilde{\rho}' \otimes \tilde{\alpha}'
\left( \frac{\partial s_{v'}}{\partial x_{u-1}} \right)
\end{array}
\right) .
\]
By (\ref{r-and-s}), we have
\[
M_G = A M_{G'}
\]
where $A = \left(\rho'(A_{i,k})\right)$
is a $n v \times n v'$ matrix.
For $I = (i_1,i_2,\ldots,i_{n(u-1)})$,
as is easily shown,
\[
\det M_{G}^I =
\det \left( A^I M_{G'} \right) =
\sum_K
\pm \left( \det A_K^I \right)
\left( \det M_{G'}^K \right)
\]
where $K = (k_1,k_2,\ldots,k_{n(u-1)})$ and
$A_K^I$ is the matrix consisting of
the $k_1, k_2,$ $\ldots,$ 
$k_{n(u-1)}$-th columns of $A^I$.
It follows that
if $\det M_{G'}^I$ has a common divisor $P$
for all $I$,
then so does $\det M_{G}^I$.
Moreover, the denominator of $\Delta_{G,\rho}$
is equal to that of $\Delta_{G',\rho'}$.
This completes the proof.
\end{proof}

The corresponding fact about the Alexander polynomial 
is well known. 
Let $G(K)$ be the knot group $\pi_1(S^3 - K)$ 
of a knot $K$ in $S^3$. 
For any knots $K,K'$, 
if there exists a surjective homomorphism 
from $G(K)$ to $G(K')$,
then the Alexander polynomial of $K$ is divisible by 
that of $K'$. 
Murasugi mentions that 
if there exists a surjective homomorphism 
from a knot group $G(K)$ to the trefoil knot group, 
then the twisted Alexander polynomial 
of $G(K)$ is divisible by 
that of the trefoil knot group. 
The main theorem is a generalization of the above. 

We will now make a few remarks about geometric settings 
in which surjective homomorphisms arise. 
First we consider the case of degree one maps. 
Let $X$ and $Y$ be $d$-dimensional compact manifolds. 
Suppose that $f \co X\to Y$ is a degree one map. 
It is easy to see that 
its induced homomorphism $f_* \co \pi_1(X)\to \pi_1(Y)$ 
is a surjective homomorphism.

In the knot group case, 
there exist the following situations except for degree 1 maps. 
First, 
there exists a surjective homomorphism 
from any knot group to the trivial knot group 
which is the infinite cyclic group. 
Secondly, 
if a knot $K$ is a connected sum of $K_1$ and $K_2$, 
then its knot group $G(K)$ is an amalgamated product 
of $G(K_1)$ and $G(K_2)$. 
Then there exists a surjection 
from $G(K)$ to each factor group. 
Thirdly, if a knot $K$ is a periodic knot of order $n$, 
then there exists a surjective homomorphism 
from $G(K)$ to $G(K_*)$ 
where $K_*$ is its quotient knot of $K$. 

%%%%%%%%%%%%%%%%%%%%%%%%%%%%%%%
\section{Another proof from the view point of 
the Reidemeister torsion}\label{geo-proof} 

In this section,
we prove our theorem in the knot group case. 
It is done by using the Mayer-Vietoris argument
of the Reidemeister torsion. 

Here let us consider a knot $K$ in $S^3$ 
and its exterior $E(K)$. 
For the knot group $G(K)=\pi_1E(K)$, 
we choose and fix a Wirtinger presentation
$$
G(K)=
\langle x_1,\ldots,x_u~|~r_1,\ldots,r_{u-1}\rangle.
$$
The abelianization homomorphism 
$$
\alpha_K \co G(K)\to H_1(E(K),\Bbb Z)\cong\Bbb Z
=\langle t\rangle
$$
is given by
$\alpha_K(x_1)=\cdots=\alpha_K(x_u)=t$.
If we have no confusion, 
we write simply $\alpha$ for $\alpha_K$ 
as in the previous section. 
In this section, 
we take a unimodular representation 
$\rho \co G(K) \to SL(n; {\mathbb F})$ over a field ${\mathbb F}$. 
As in the definition of the twisted Alexander polynomial, 
we consider the tensor representation 
$$
\rho\otimes\alpha \co 
G \to GL(n;\F[t,t^{-1}])\subset GL(n;\F(t)).
$$
Here 
$\F(t)$ denotes the rational function field
over $\F$. 
If $\rho\otimes\alpha$ is an acyclic representation over $\F(t)$, 
that is, 
all homology groups over $\F(t)$ 
of $E(K)$ twisted by $\rho\otimes\alpha$ 
are vanishing, 
then the Reidemeister torsion of $E(K)$ 
for $\rho\otimes\alpha$ can be defined. 
Furthermore the following equality holds. 
See \cite{KL, Kitano} for more details of definitions and proofs. 

\begin{theorem}
If $\rho\otimes\alpha$ is an acyclic representation,
then
we have
$$
\tau_{\rho\otimes\alpha}(E(K))
=\Delta_{G(K),\rho}(t)
$$
up to a factor
$\pm t^{nk}~(k\in\Bbb Z)$ if $n$ is odd, 
and up to only $t^{nk}$ if $n$ is even.
\end{theorem}

From this theorem, 
we prove the main theorem as divisibility of 
the Reidemeister torsion 
in the knot group case. 
Here we take a surjective homomorphism
$\varphi \co G(K)\to G(K')$.
By changing the orientation of meridians if we need,
we may assume that $\alpha_{K'}\circ \varphi=\alpha_K$.
Let $\rho' \co G(K')\to SL(n;\F)$ be a representation.
For simplicity,
we write the composition $\rho = \rho'\circ\varphi$. 

Now we consider 2-dimensional CW-complexes
$X(K)$ and $X(K')$
defined by their Wirtinger presentations.
It is well-known that
these complexes are simple homotopy equivalent
to the knot exteriors.
Then these Reidemeister torsions of $X(K)$ and $X(K')$
are equal to the twisted Alexander polynomials respectively. 
Here we consider twisted homologies of these complexes
by using their CW-complex structure.
The coefficient $V$ is a $2n$-dimensional vector space 
over a rational function field $\F(t)$.
When $V$ is regarded as a $G(K)$-module by using $\rho$, 
it is denoted by $V_\rho$. 

The 
surjective 
homomorphism $\varphi$ induces a chain map
$\varphi_* \co C_*(X(K),V_{\rho})\to
C_*(X(K'),V_{\rho'})$.
We take a tensor representation 
$\rho\otimes\alpha_K \co G(K)\to GL(n;\F(t))$. 
Assume that
$\rho\otimes \alpha_K$ and $\rho'\otimes \alpha_{K'}$
are acyclic representations.
Then we can prove the following.

\begin{theorem}\label{geo-main}
The quotient 
$\tau(X(K);V_{\rho\otimes\alpha_K})/\tau(X(K');
V_{\rho'\otimes\alpha_{K'}})$
is a polynomial in $\F[t,t^{-1}]$.
\end{theorem}

We show the following proposition first. 

\begin{proposition}
The chain map
\[
\varphi_* \co C_*(X(K),V_{\rho\otimes\alpha_{K}})
\to C_*(X(K'),V_{\rho'\otimes\alpha_{K'}})
\]
is surjective.
\end{proposition}

\begin{proof}
It is clear that $\varphi$ induces
an isomorphism on the 0-chains,
and a surjection on the 1-chains.
Then we only need to prove this proposition on the
2-chains.

We take a non-trivial 2-chain $z\in 
C_2(X(K'),V_{\rho'\otimes\alpha_{K'}})$.
By the acyclicity of 
the chain complex 
$C_*(X(K'),V_{\rho'\otimes\alpha_{K'}})$,
the boundary map
$\partial \co C_2(X(K'),V_{\rho'\otimes\alpha_{K'}})$ 
$\to C_1(X(K'),V_{\rho'\otimes\alpha_{K'}})$
is injective.
Then the image $\partial z$ is non-trivial in $C_1$. 
On the other hand,
by the surjectivity of
$$
\varphi \co C_1(X(K),V_{\rho\otimes\alpha_{K}})
\to
C_1(X(K'),V_{\rho'\otimes\alpha_{K'}}),
$$
there exists a 2-chain $w\in
C_2(X(K),V_{\rho\otimes\alpha_{K}})$
such that $\varphi_*(w)=z$.
By the commutativity of maps, 
in $C_2$ 
$$
\varphi_*(\partial w)=\partial\varphi_*(w)=\partial \partial
z=0.
$$
Then we have $\partial w=0$.
%Hence b
By the acyclicity, 
there exists $\tilde w\in
C_*(X(K),V_{\rho\otimes\alpha_K})$
such that $\partial\tilde w=w$.
Again by the commutativity,
$\varphi(\tilde w)=z$.
Therefore $\varphi_*$ is surjective.
\end{proof}

\begin{proof}[Proof of Theorem \ref{geo-main}]
From the above proposition, 
we can take the kernel $D_*$
of this chain map $\varphi_*$
and obtain a short exact sequence
$$0\to D_* \to C_*(X(K),V_{\rho\otimes\alpha_K})
\to C_*(X(K'),V_{\rho'\otimes\alpha_{K'}})
\to 0.$$

Here we recall the following fact. 
For a short exact sequence $0\to C'_*\to C_*\to C''_* \to 0$ 
of finite chain complexes, 
if two of them are acyclic complexes, 
then the third one is also acyclic. 
Furthermore, the torsion satisfies 
$$
\tau(C_*)=\tau(C'_*)\tau(C''_*)
$$ 
up to some factor. 

By applying the property of the product of torsion, 
we have
$$\tau(X(K);V_{\rho\otimes\alpha_K})
=\tau(X(K');V_{\rho'\otimes\alpha_{K'}})
\tau(D;V_{\rho\otimes\alpha_K}).$$

We only need to prove that
$\tau(D;V_{\rho\otimes\alpha_K})$ is a polynomial.
From the definition we see that $D_0$ vanishes, 
since 
$$
\varphi_* \co C_0(X(K),V_{\rho\otimes\alpha_K})
\to C_0(X(K'),V_{\rho'\otimes\alpha_{K'}})
$$ is isomorphism.
Hence by definition,
its torsion is the determinant of
$D_2\to D_1$.
Therefore it is a polynomial.
\end{proof}

\begin{remark}
By a similar argument, 
we can prove that 
if $\varphi \co G(K)\to G(K')$ is an injective homomorphism, 
then 
$\tau(X(K');V_{\rho\otimes\alpha_{K'}})
/\tau(X(K);V_{\rho\otimes\alpha_K})$
is a polynomial.
\end{remark}

%%%%%%%%%%%%%%%%%%%%%%%%%%%%%%%%%%%%%%%%%%%%%%%%%%%
\section{Examples}
In this section, we show some examples 
of the twisted Alexander polynomials and 
an application of Theorem \ref{main-theorem}. 
We consider the problem:\ Is 
there a surjective homomorphism 
from $G(K)$ to $G(K')$ for two given knots $K,K'$? 
The problem has been investigated by Murasugi 
when $K'$ is the trefoil knot $3_1$ (c.f.\ \cite{mur}). 
Here we study the problem 
in case when $K'$ is the figure eight knot $4_1$. 
The numbering of the knots follows
that of Rolfsen's book \cite{R}.

If the classical Alexander polynomial of $K$
can not be divided by that of $K'$,
we 
know that there are no surjective homomorphisms 
from $G(K)$ to $G(K')$.
In the knot table in \cite{R}, 
up to $9$ crossings, 
the classical Alexander polynomial of 
each knot is not divisible by that of $G(4_1)$ 
except for
$8_{18}, 8_{21}, 9_{12}, 9_{24}, 9_{37},$ 
$9_{39}$ and $9_{40}$. 
That is to say, 
except for
$8_{18}, 8_{21}, 9_{12}, 9_{24}, 9_{37},$ 
$9_{39}$ and $9_{40}$, 
there exists no surjective homomorphisms from
such a knot group to $G(4_1)$.

Next, we consider a representation 
$\rho \co G(K) \to SL(2;{\mathbb Z} / p {\mathbb Z})$ 
and the twisted Alexander polynomial
associated to $\rho$. 
Theorem \ref{main-theorem} says that
if the numerator of $\Delta_{G(K),\rho}$ 
for all representations
$\rho \co G(K) \to SL(2;{\mathbb Z} / p {\mathbb Z})$
for some fixed prime $p$ cannot be divided
by the numerator of $\Delta_{G(K'),\rho'}$
for a certain representation
$\rho' \co G(K') \to SL(2;{\mathbb Z} / p {\mathbb Z})$,
then there exists no surjective homomorphisms 
from $G(K)$ to $G(K')$.

Let us compute the twisted Alexander polynomials
$\Delta_{G(4_1),\rho'}$ for a certain representation
$\rho' \co G(4_{1}) \to SL(2;{\mathbb Z} / 7 {\mathbb Z})$.
The knot group $G(4_{1})$ admits a presentation
\[
G(4_{1}) = \langle
x_1,x_2,x_3,x_4 \,\, | \,\,
x_4 x_2 x_4^{-1} x_1^{-1},\,
x_1 x_2 x_1^{-1} x_3^{-1},\,
x_2 x_4 x_2^{-1} x_3^{-1}
\rangle .
\]
We can check easily
that the following is a representation of $G(4_{1})$:
\[
\rho' (x_1) =
\left(
\begin{array}{cc}
1 & 1 \\
0 & 1
\end{array}
\right), \,
\rho' (x_2) =
\left(
\begin{array}{cc}
1 & 0 \\
3 & 1
\end{array}
\right), \,
\]\[
\rho' (x_3) =
\left(
\begin{array}{cc}
4 & 4 \\
3 & 5
\end{array}
\right), \,
\rho' (x_4) =
\left(
\begin{array}{cc}
2 & 4 \\
5 & 0
\end{array}
\right).
\]
Then we obtain the Alexander matrix:
\[
M =
\left(
\begin{array}{cccccccc}
6 & 0 & 2 t & 4 t & 0 & 0 & 6 t + 1 & 6 t \\
0 & 6 & 5 t & 0 & 0 & 0 & 0 & 6 t + 1 \\
3 t + 1 & 3 t & t & t & 6 & 0 & 0 & 0 \\
4 t & 2 t + 1 & 0 & t & 0 & 6 & 0 & 0 \\
0 & 0 & 3 t + 1 & 3 t & 6 & 0 & t & 0 \\
0 & 0 & 4 t & 2 t + 1 & 0 & 6 & 3 t & t 
\end{array}
\right) 
\]
The numerator $P$ of
the twisted Alexander polynomial $\Delta_{G(4_1),\rho'}$
is the determinant of $M_4$
obtained from $M$ by removing the last two columns.
Then we get
\[
P = t^4 + t^3 + 3 t^2 + t +1 .
\]
Moreover,
we calculate the numerator of
the twisted Alexander polynomials of $G(8_{21})$
for all representations
$G(8_{21}) \to SL(2;{\mathbb Z} / 7 {\mathbb Z})$
and get $24$ polynomials. 
These calculations are made by author's computer program 
and the same results are obtained 
by Kodama Knot program \cite{Kodama}. 
None of them can be divided by $P$,
so we conclude that
there exists no surjective homomorphisms 
from $G(8_{21})$ to $G(4_1)$.
By similar arguments 
using $SL(2;{\mathbb Z} / p {\mathbb Z})$-representations 
for $p=5,7$, 
we get the conclusion that 
there exists no surjective homomorphisms 
from $G(9_{12}),G(9_{24}),G(9_{39})$ to $G(4_1)$.
On the other hand,
$8_{18}$ is a periodic knot of order $2$ 
with quotient knot $4_1$. 
Furthermore, $G(9_{37})$ has a presentation 
\begin{eqnarray*}
\lefteqn{G(9_{37}) = }\\
&& \left\langle
\begin{array}{c}
y_1,y_2,y_3,y_4,y_5, \\
y_6,y_7,y_8,y_9 ~
\end{array}
\left| ~ 
\begin{array}{c}
y_8 y_1 y_8^{-1} y_2^{-1}, 
y_7 y_2 y_7^{-1} y_3^{-1}, 
y_9 y_4 y_9^{-1} y_3^{-1}, 
y_3 y_4 y_3^{-1} y_5^{-1}, \\
y_1 y_6 y_1^{-1} y_5^{-1}, 
y_5 y_6 y_5^{-1} y_7^{-1}, 
y_2 y_7 y_2^{-1} y_8^{-1}, 
y_4 y_9 y_4^{-1} y_8^{-1} 
\end{array}
\right.
\right\rangle 
\end{eqnarray*}
and the following mapping 
$\varphi \co  G(9_{37}) \to G(4_1)$ 
is a surjective homomorphism: 
\[
\begin{array}{c}
\varphi(y_1) = x_2,
\varphi(y_2) = x_3,
\varphi(y_3) = x_1 x_4 x_1^{-1},
\varphi(y_4) = x_3,
\varphi(y_5) = x_1, \\
\varphi(y_6) = x_1^{-1} x_4 x_1,
\varphi(y_7) = x_4,
\varphi(y_8) = x_1,
\varphi(y_9) = x_4.
\end{array}
\]
Similarly, we can give an explicit surjective homomorphism from
the knot group $G(9_{40})$ to $G(4_1)$.  Thus we have surjective
homomorphisms from knot groups $G(8_{18})$, $G(9_{37})$, $G(9_{40})$
to $G(4_1)$.  Hence we can determine whether or not there exists a
surjective homomorphism from the group of each knot with up to $9$
crossings to $G(4_1)$.

In \cite{KS}, 
we see a complete list of 
whether there exists a surjective homomorphism 
between knot groups for $10$ crossings and less. 

\medskip
{\bf Acknowledgements}\qua
The first author is supported 
in part by Grand-in-Aid for Scientific Research
(No.\ 14740037),
The Ministry of Education, Culture,
Sports, Science and Technology, Japan.
The second author is supported by
the 21 century COE program
at Graduate School of Mathematical Sciences,
the University of Tokyo. 
The authors would like to express their thanks
to Prof.\ Sadayoshi Kojima and Prof.\ Dieter Kotschick 
for their useful comments.

%%%%%%%%%%%%%%%%%%%%%%%%%%%%%%%%%%%%%%%%%%%%%%%%%%%%%%%%

%\bibliographystyle{amsplain}
\bibliographystyle{gtart}

\Addresses\recd

\end{document}